\documentclass[10pt]{article}
\usepackage{amssymb,amsmath,latexsym}

\newtheorem{thm}{Theorem}[section]

\newtheorem{remark}[thm]{Remark}
\newtheorem{defn}[thm]{Definition}

\newtheorem{example}[thm]{Example}
\newtheorem{lemma}[thm]{Lemma}

\numberwithin{equation}{section}

\def\pf{\noindent{\bf Proof.} }
\def\qed{{\hfill $\Box$ \bigskip}}

\def\topdot{{\leavevmode
      \raise.2ex\hbox{${\cdot}$}}}

\def\EE{{\mathcal E}}
\def\FF{{\mathcal F}}

\def\P{{\mbox{\bf P}}}
\def\E{{\mbox{\bf E}}}
\def\1{{\mbox{\bf 1}}}

\def\loc{{\rm loc}}

\def\<{\langle}
\def\>{\rangle}

\newcommand{\Rn}{\mathbb{R}^n}

\newcommand{\R}{\mathbb{R}}

\makeatletter

\@addtoreset{equation}{section} \makeatother

\begin{document}
\title{\bf On notions of harmonicity }
\author{
{\bf Zhen-Qing  Chen}\footnote{The research of this author is
supported in part by NSF Grant DMS-0600206.}}
\bigskip
\date{ }
\maketitle

\begin{abstract}
In this  paper, we address the equivalence of the analytic and
probabilistic notions of harmonicity in the context of general
symmetric Hunt processes on locally compact separable metric spaces.
Extensions to general symmetric right processes on Lusin spaces
including infinite dimensional spaces are  mentioned at the end of
this paper.
 \end{abstract}



\bigskip

\begin{doublespace}

\section{Introduction}

It is known that a function $u$ being harmonic in a domain $D\subset
\Rn$ can be defined or characterized by $\Delta u=0$ in $D$ in the
distributional sense, that is, $u\in W^{1,2}_\loc(D):=\left\{v\in
L^2_\loc (D) \mid \nabla v \in L^2_\loc (D)\right\}$ so that
$$ \int_{\Rn} \nabla u (x) \cdot \nabla v(x) dx =0
\qquad \hbox{for every } v\in C^\infty_c(D). $$
It is equivalent to the following averaging property by running a Brownian motion
$X$: for every relatively compact subset $U$ of $D$,
$$ u(X_{\tau_U})\in L^1(\P_x) \qquad \hbox{and} \qquad
u(x) = \E_x \left[ u(X_{\tau_U})\right] \quad \hbox{for every } x\in
U.
$$
Here $\tau_U:=\inf\left\{t\geq 0: X_t \notin U\right\}$. Recently
there are interests (e.g. \cite{BKK}) resulting from several areas
of mathematics in knowing whether the above two notions of
harmonicity remain equivalent in a more general context, such as for
diffusions on fractals (see \cite{BBKT}) and for discontinuous
processes including symmetric L\'evy processes. For instance, due to
their importance in theory and in applications, there has been
intense interest recently in studying discontinuous processes and
non-local (or integro-differential) operators, by
 both analytical and probabilistic approaches.
 See, e.g.,  \cite{CK, CK2} and the references therein. So it is important
 to identify the connection between the analytic and probabilistic notions
 of harmonic functions.

\medskip

In this  paper, we address the question of the equivalence of the
analytic and probabilistic notions of harmonicity in the context of
symmetric Hunt processes on local compact separable metric spaces.
Let $X$ be an $m$-symmetric Hunt process on a locally compact
separable metric space $E$ whose associated Dirichlet form $(\EE,
\FF)$ is regular on $L^2(E; m)$. Let $D$ be an open subset of $E$
and $\tau_D$ is the first exit time from $D$ by $X$. Motivated by
the example at the beginning of this section, loosely speaking (see
next section for precise statements), there are two ways to define a
function $u$ being harmonic in $D$ with respect to $X$: (a)
(probabilistically) $t\mapsto u(X_{t\wedge \tau_D})$ is a
$\P_x$-uniformly integrable martingale for quasi-every $x\in D$; (b)
(analytically) $\EE (u, g)=0$ for  $g\in \FF\cap C_c(D)$.
 We will show in Theorem \ref{T:7} below that these two definitions
 are equivalent.
Note that even in the Brownian motion case,  a function $u$ that is
harmonic in $D$ is typically not in the domain $\FF$ of the
Dirichlet form. Denote by $\FF^D_\loc$   the family of functions $u$
on $E$  such that for every relatively compact open subset $D_1$ of
$D$, there is a function $f\in \FF$ so that $u=f$ $m$-a.e. on $D_1$.
To show these two definitions are equivalent, the crux of the
difficulty is to
\begin{description}
\item{(i)} appropriately extend the definition of $\EE(u, v)$ to functions $u$ in
$ \FF^D_\loc$ that satisfy some minimal integrability condition when
$X$ is discontinuous  so that $\EE(u, v)$ is well defined for every
$v\in \FF\cap C_c(D)$;

\item{(ii)} show that if $u$ is harmonic in $D$ in the probabilistic sense,
then $u\in \FF^D_\loc$ and $\EE(u, v)=0$ for every $v\in \FF\cap
C_c(D)$.
\end{description}
If one assumes a priori that $u\in \FF$, then the equivalence of (a)
and (b) is  easy to establish. See Remarks \ref{R:5}(i) and
\ref{R:2.7} below.

In next section, we give precise definitions, statements of the main
results and their proofs. Three examples are given to illustrate the
main results of this paper. Extensions to general symmetric right
processes on Lusin spaces including infinite dimensional spaces are
mentioned at the end of this paper. We use ``:=" as a way of
definition. For two real numbers $a$ and $b$, $a\wedge b :=\min\{a,
b\}$.

\section{Main results}

   Let $X=(\Omega,\FF_{\infty},\FF_t,
X_t,\zeta,\P_x, x\in E)$ be an $m$-symmetric Hunt   process on a
locally compact separable metric space $E$, where $m$ is a positive
Radon measure on $E$ with full topological support. A cemetery state
$\partial$ is added to $E$ to form $E_\partial:=E\cup\{\partial\}$
as its one-point compactification and  $\Omega$ is the totality of
right-continuous, left-limited
 sample paths from $[0,\infty[$ to
$E_\partial$ that hold the value $\partial$ once attaining it.  For
any $\omega\in\Omega$, we set $X_t(\omega):=\omega(t)$.  Let
$\zeta(\omega):=\inf\{t\geq0\,\mid\, X_t(\omega)=\partial\}$  be the
life time of $X$. As usual, $\FF_{\infty}$ and $\FF_t$ are the
minimal augmented $\sigma$-algebras obtained from
$\FF_{\infty}^0:=\sigma\{X_s\,\mid\, 0\leq s<\infty\}$ and
$\FF_t^0:=\sigma\{X_s\,\mid\, 0\leq s\leq t\}$ under $\{\P_x: x\in
E\}$. For a Borel subset $B$ of $E$, $\tau_B:=\inf\{t>0 \mid
X_t\notin B\}$ (the {\it exit time} of $B$) and
$\sigma_B:=\inf\{t\geq 0 \mid X_t\in B\}$ (the {\it entrance time} of $B$) are $(\FF_t)$-stopping times.

The transition semigroup $\{P_t: t\ge 0\}$ of $X$ is defined by
$$
P_tf(x):=\E_x[f(X_t)]=\E_x[f(X_t): t< \zeta],\qquad t\ge 0.
$$
Each  $P_t$ may be viewed as an operator on $L^2(E, m)$, and taken as a whole
these operators form a strongly continuous semigroup of
self-adjoint contractions. The Dirichlet form associated with $X$ is the
bilinear form
\begin{equation}\label{e:1.1}
{\cal E}(u, v):=\lim_{t\downarrow 0}t^{-1}(u-P_tu, v)_m
\end{equation}
defined on the space
\begin{equation}\label{e:1.2}
{\cal F}:=\left\{u\in L^2(E; m)\,\Big|
\,\sup_{t>0}\,\,t^{-1}(u-P_tu, u)_m<\infty \right\}.
\end{equation}
Here we use the notation $(f,g)_m:=\int_E f(x)g(x)\, m(dx)$. We
assume that $(\EE, \FF)$ is a regular Dirichlet form on $L^2(E; m)$;
that is, $C_c(E)\cap \FF$ is dense both in $(C_c(E), \| \cdot
\|_\infty)$ and in $(\FF, \EE_1)$. Here $C_c(E)$ is the space of
continuous functions with compact support in $E$ and $\EE_1(u,
u):=\EE (u, u)+(u, u)_m$. However  to ensure a  wide scope of
applicability, we do {\it not} assume that  the process $X$ (or
equivalently,  its associated  Dirichlet form $(\EE, \FF)$) is
$m$-irreducible.

We refer readers to \cite{CF} and \cite{FOT}
for the following known facts. The
extended Dirichlet space $\FF_e$ is the space of all functions
$f$ on $E$ so that there is an $\EE$-Cauchy sequences $\{f_n, n\geq
1\}\subset \FF$ so that $f_n$ converges to $f$ $m$-a.e. on $E$. For
such an $f\in \FF_e$, $\EE(f, f):=\lim_{n\to \infty} \EE (f_n,
f_n)$. Every $f\in \FF_e$ admits a quasi-continuous version (cf.
\cite[Theorem 2.1.7]{FOT}). Throughout this paper, we always   assume
that every function
 in $\FF_e$ is represented by its quasi-continuous version, which is unique
 up to a set of zero capacity (that is, quasi-everywhere, or q.e. in abbreviation).
 We adopt  the convention that
 any function $f$ defined on $E$ is extended to $E_\partial$ by taking $f(\partial )=0$
 and that  $X_\infty (\omega):=\partial$ for every $\omega \in \Omega$.
 It is known that $\FF_e\cap L^2(E; m] =\FF$. The extended Dirichlet
form $(\EE, \FF_e)$ admits the following Beurling-Deny decomposition
(cf. \cite[Theorem 4.3.3]{CF} or \cite[Theorem 5.3.1]{FOT}):
$$ \EE (u, u)=\EE^{(c)}(u, u) + \frac12 \int_{E\times E} (u(x)-u(y))^2 J(dx, dy)
+ \int_E u(x)^2 \kappa (dx),
$$
 where $\EE^{(c)}$ is the strongly local part of $(\EE, \FF)$,  $J$  the jumping measure
 and $\kappa$ the killing measure of $(\EE, \FF)$ (or, of $X$).
 For $u, v\in \FF_e$, $\EE^{(c)}(u, v)$ can  also be expressed
 by the mutual energy measure
 $\frac12 \mu^c_{\<u, v\>}(E)$, which is the signed Revuz measure associated with
 $\frac12 \< M^{u, c}, M^{v, c}\>$. Here for $u\in \FF_e$,
 $M^{u, c}$ denotes the continuous martingale
 part of the square integrable martingale additive functional $M^u$ of $X$ in the Fukushima's
 decomposition (cf. \cite[Theorem 5.2.2]{FOT}) of
 $$u(X_t)-u(X_0)=M^u_t+N^u_t, \qquad t\geq 0,
 $$
 where $N^u$ is continuous additive functional of $X$ having zero energy.
 When $u=v$, it
 is  customary to write  $\mu^c_{\<u, u\>}$  as
 $\mu^c_{\<u\>}$. The measure $\mu^c_{\<u, v\>}$ enjoys the strong local property in
 the sense that if $u\in \FF_e$ is constant on a nearly Borel quasi-open set $D$,
 then $\mu^c_{\<u, v \>}(D)=0$ for every $v\in \FF_e$
 (see \cite[Proposition 4.3.1]{CF}).
 For $u\in \FF$, let $\mu_{\<u\>}$ be the Revuz measure of $\<M^u\>$. Then
 it holds that
 $$ \EE(u, u)=\frac12 \mu_{\<u\>}(E)+\frac 12 \int_E u(x)^2 \kappa (dx).
 $$
For an open subset $D$ of $E$, we use $X^D$ to denote the
subprocess of $X$ killed upon leaving $D$. The Dirichlet form of
$X^D$ on $L^2(D; m)$ is $(\EE, \FF^D)$, where $\FF^D:=\{u\in \FF \mid
u=0 \hbox{ q.e. on } D^c\}$. It is known (cf. \cite[Theorem 3.3.9]{CF} or \cite[Theorem
4.4.3]{FOT} that $(\EE, \FF^D)$ is a regular Dirichlet form on
$L^2(D; m)$. Let $\FF^D_e:=\{u\in \FF_e \mid u=0 \hbox{ q.e. on }
D^c\}$. Then $\FF^D_e$ is the extended Dirichlet space of $(\EE,
\FF^D)$ (see Theorem 3.4.9 of \cite{CF}).
 A function $f$ is said to be locally in $\FF^D$, denoted as $f\in \FF^D_{\loc}$,
 if for every relatively compact subset $U$ of $D$, there is a function $g\in \FF^D$
 such that $f=g $ $m$-a.e. on $U$.  Every $f\in \FF^D_\loc$ admits an $m$-version
 that is quasi-continuous on $D$.
 Throughout this paper, we always   assume
that every function
 in $\FF^D_\loc$, when restricted to $D$ is represented by its quasi-continuous  version.
 By the strong local property of $\mu^c_{\<u, v\>}$
 for $u, v \in \FF$,
 $\mu^c_{\<u, v\>}$ is well defined on $D$ for every $u, v\in \FF^D_\loc$.
 We use $L^\infty_\loc (D; m)$ to
 denote the $m$-equivalent class of locally bounded functions on $D$.

 Let   $(N(x, dy), H)$ be a L\'evy system of   $X$ (cf. \cite{BJ} or \cite{FOT}).
  Then
  $$J(dx, dy)=N(x, dy) \mu_H(dx) \qquad \hbox{and} \qquad
  \kappa (dx):=N(x, \partial ) \mu_H(dx),
  $$
  where $\mu_H$ is the Revuz measure of the positive continuous additive functional
   $H$ of $X$.

\begin{defn}\label{D:1}
\rm  Let $D$ be an open subset of $E$. We say a function
$u$ is {\it harmonic} in $D$ (with respect to the process $X$) if
  for every relatively compact open subset $U$ of $D$,
  $t\mapsto u(X_{t\wedge \tau_U})$
 is a uniformly integrable $\P_x$-martingale for q.e. $x\in U$.
\end{defn}

\medskip

To derive an analytic characterization
of harmonic functions in $D$ in terms of an extension of  quadratic form $(\EE, \FF)$,
 we need some preparation.
Let $r_t$ denote the time-reversal operator defined on the path space
$\Omega$  of $X$ as follows: For  $\omega\in\{t<\zeta\}$,
$$
r_t(\omega)(s)= \begin{cases} \omega((t-s){-})  & \hbox{if } 0\le s<
t,\\ \omega(0)  & \hbox{if } s\ge t.
\end{cases}
$$
(It should be borne in mind that the restriction of the measure
$\P_m$ to $\FF_t$ is invariant under $r_t$ on $\Omega \cap \{ \zeta
>t\}$.)

 \begin{lemma}\label{L:2} If $u\in \FF_e$ has $\EE(u, u)=0$, then
 $$\P_x \left( u(X_t)=u(X_0)
  \hbox{ for every } t\geq 0\right)=1 \qquad \hbox{for q.e. } x\in E.
  $$
  In other words, for q.e. $x\in E$, $E_x:=\{y\in E: u(y)=u(x)\}$ is an invariant
  set with respect to the process $X$ in the sense that $\P_x(X [0, \infty)\subset E_x)=1$. This in particular implies that, if, in addition, $\P_x(\zeta<\infty)>0$ for q.e.
  $x\in E$, then $u=0$ q.e. on $E$.
 \end{lemma}

 \pf It is known (see, e.g., \cite[Theorem 6.6.2]{CF}) that the following Lyons-Zheng's
 forward-backward martingale decomposition holds for $u\in \FF_e$:
  $$ u(X_t)-u(X_0)= \frac12 M^u_t -\frac12 M^u_t \circ r_t \quad \P_m \hbox{-a.e. on }
  \{t<\zeta \}.
  $$
 As $\mu_{\<u\>}(E)\leq 2 \EE(u, u)=0$, we have $M^u=0$ and so
 $u(X_t)=u(X_0)$ $\P_m$-a.s. on $\{t<\zeta \}$ for every $t>0$.
 This implies via Fukushima's decomposition that
 $N^u=0$ on $[0, \zeta)$ and hence on $[0, \infty)$ $\P_m$-a.s. Consequently,
 $ \P_x \left( u(X_t)-u(X_0)=M^u_t+N^u_t=0  \hbox{ for every } t\geq
 0\right)=1$ for q.e. $x\in E$.
This proves the lemma. \qed

 \medskip

Since $(\EE, \FF)$ is a regular Dirichlet form on $L^2(E; m)$, for
  any relatively compact open sets $ U, V$ with
$\overline U\subset V$, there is $\phi \in \FF\cap C_c(E)$ so that
$\phi =1$ on $U$ and $\phi =0$ on $V^c$. Consequently,
\begin{equation}\label{e:J1}
J(U, V^c) = \int_{U\times V^c} (\phi (x)-\phi (y))^2 J(dx, dy)
\leq 2 \EE (\phi, \phi)<\infty.
\end{equation}

 For an open set $D\subset E$,
 consider the following two conditions for function $u$ on $E$.
For any relatively compact open sets $U, V$ with
 $\overline U \subset V \subset \overline V \subset D$,
 \begin{equation}\label{e:cond1}
 \int_{U\times  (E\setminus V)} |u(y)| J(dx, dy) <\infty
 \end{equation}
 and
\begin{equation}\label{e:cond2}
\1_U (x) \E_x \left[ \big((1-\phi_V ) |u|\big)(X_{\tau_U})\right]
   \in \FF^U_e,
 \end{equation}
 where $\phi_V\in C_c(D)\cap \FF$ with $0\leq \phi_V\leq 1$ and $\phi_V=1$
 on $V$. Note that both conditions \eqref{e:cond1} and \eqref{e:cond2} are automatically
 satisfied when $X$ is a diffusion since in this case
 the jumping measure $J$
 vanishes and $X_{\tau_U}\in \partial U $ on $\{\tau_U<\zeta\}$.
In view of \eqref{e:J1}, every bounded function $u$ satisfies the
  condition \eqref{e:cond1}. In fact by the following lemma, every bounded function $u$
  also satisfies the condition \eqref{e:cond2}.

  \begin{lemma}\label{L:2.3} Suppose that $u$ is a function on $E$
  satisfying condition \eqref{e:cond1} and that
  for any relatively compact open sets $U, V$ with
 $\overline U \subset V \subset \overline V \subset D$,
  \begin{equation}\label{e:cond3}
  \sup_{x\in U} \E_x \left[ \big(\1_{V^c} |u|\big)
   (X_{\tau_U})\right]<\infty.
  \end{equation}
 Then   \eqref{e:cond2} holds for $u$.
\end{lemma}

In many concrete cases such as in Examples \ref{E:8}-\ref{E:10}
below, one can show that condition \eqref{e:cond1} implies condition
\eqref{e:cond3}. To prove the above lemma, we need the following
result.
Observe that the process $X$ is not assumed to be
transient.

\begin{lemma}\label{L:2.4}
 Suppose that $\nu$ is a smooth measure on $E$, whose corresponding
positive continuous additive functional (PCAF) of $X$ is denoted as
$A^\nu$. Define $G\nu (x):= \E_x[A^\nu_\zeta]$. If $\int_E G\nu (x)
\nu (dx)<\infty$, then $G\nu \in \FF_e$. Moreover,
\begin{equation}\label{e:1}
\EE(G \nu, u)=\int_E u(x) \nu(dx) \qquad \hbox{for every } u\in
\FF_e.
\end{equation}
\end{lemma}

\pf  First assume that $m(E)<\infty$. It is easy to check directly
that
 $\{x\in E: \E_x [ A^\nu_\zeta]>j\}$
is finely open for every integer $j\geq 1$. So $K_j:=\{G\nu \leq
j\}$ is finely closed. Since $G\nu <\infty$ $\nu$-a.e. on $E$, we
have $\nu (E\setminus \cup_{j=1}^\infty K_j)=0$. Define
$\nu_j:=\1_{K_j} \nu$. Clearly for $x\in K_j$, $G\nu_j(x) \leq G\nu
(x)\leq j$, while for $x\in K_j^c$,
$$G\nu_j (x)=\E_x \left[ \int_0^\zeta
\1_{K_j}(X_s) dA^\nu_s\right]= \E_x \left[ G\nu_j
(X_{\sigma_{K_j}})\right]\leq j.
$$
So $f_j:=G\nu_j\leq j$ on $E$ and hence is in $L^2(E; m)$. Since
by \cite[Theorem 4.1.1]{CF} or \cite[Theorem 5.1.3]{FOT}
\begin{equation}\label{e:2}
 \lim_{t\to 0} \frac1t (f_j-P_tf_j, \, f_j)_m
  = \lim_{t\to 0} \frac 1t \E_{f_j\cdot m} \left[ A^{\nu_j}_t \right]
  = \int_E f_j (x) \nu_j(dx) \leq \int_E G\nu (x) \nu (dx)<\infty,
\end{equation}
  we have $f_j\in \FF$ with $\EE(f_j, f_j)\leq \int_E G\nu (x) \nu (dx)$.
  The same calculation shows that for $i>j$, $f_i-f_j= \E_x \left[ A_\zeta^{\1_{K_i\setminus K_j}\cdot \nu}\right]$
  and
  $$ \EE (f_i-f_j, f_i-f_j) =  \int_{K_i\setminus K_j} (f_i-f_j) (x) \nu(dx)
  \leq \int_{K_l\setminus K_j} G\nu(x) \nu(dx),
  $$
  which tends to zero as $i, j\to \infty$; that is, $\{f_j, j\geq 1\}$ is an $\EE$-Cauchy sequence in $\FF$.  As $\lim_{j\to\infty} f_j=f$ on $E$, we conclude that $f\in \FF_e$. We deduce from
    \eqref{e:2} that
  \begin{equation}\label{e:3}
   \EE(f, f)=\lim_{j\to \infty} \EE (f_j, f_j)=\int_E G\nu (x) \nu(dx).
  \end{equation}
  Moreover, for $u\in \FF_b^+$,  by \cite[Theorem 4.1.1]{CF}
(or \cite[Theorem 5.1.3]{FOT}) and dominated convergence theorem, we have
  $$ \EE(G\nu, u)=\lim_{j\to \infty} \EE (f_j, u)
   = \lim_{j\to \infty} \lim_{t\to 0} \frac1t (f_j-P_t f_j, u) =
    \lim_{j\to \infty} \int_E u(x) \1_{K_j}(x) \nu (dx)= \int_E u(x) \nu (dx).
  $$
  Since the linear span of $\FF^+_b$ is $\EE$-dense in $\FF_e$, we have established
  \eqref{e:1}.

  For a general $\sigma$-finite measure $m$, take a strictly positive $m$-integrable
  Borel measurable function $g$ on $E$ and define $\mu=g\cdot m$. Then $\mu$ is a
  finite measure on $E$. Let $Y$ be the time-change of $X$ via measure $\mu$;
  that is, $Y_t=X_{\tau_t}$, where $\tau_t=\inf\{s>0: \int_0^s g(X_s) ds >t\}$.
 The time-changed process $Y$ is $\mu$-symmetric.
 Let $(\EE^Y, \FF^Y)$ be the Dirichlet form of $Y$ on $L^2(E; \mu)$.
   Then it is known that $\FF^Y_e=\FF_e$ and $\EE^Y=\EE$ on $\FF_e$ (see (5.2.17) of \cite{CF}).
   The measure $\nu$ is also a smooth measure with respect to process $Y$. It is easy
   to verify that the PCAF $A^{Y, \nu}$ of $Y$ corresponding to $\nu$
   is related to corresponding PACF $A^\nu$ of $X$ by
    $$ A^{Y, \nu}_t =A^\nu_{\tau_t} \qquad \hbox{for } t\geq 0 .
    $$
    In particular, we have $G^Y\nu (x)=G\nu$ on $E$. As we just proved that the lemma
     holds for $Y$, we conclude that the lemma also holds for $X$. \qed

\medskip

\noindent{\bf Proof of Lemma \ref{L:2.3}.} For relatively compact
open sets $U$, $V$ with $\overline U\subset V \subset \overline V
\subset D$ and $\phi_V\in \FF\cap C_c(D)$ with $0\leq \phi_V\leq 1$
and $\phi_V=1$ on $V$, let $f(x):=\1_U (x) \E_x \left[
\big((1-\phi_V ) |u|\big)(X_{\tau_U})\right]$, which is bounded by
condition \eqref{e:cond3}. Note that $1-\phi_V=0$ on $V$. Using
L\'evy system of $X$, we have
$$ f(x)= \E_x \left[ \int_0^{\tau_U} \left( \int_{E\setminus V} (1-\phi_V(X_s)) |u|(X_s)
N(X_s, dy)\right) dH_s  \right] \qquad \hbox{for } x\in E.
$$
Note that the Revuz measure for PCAF $t\mapsto \int_0^{t\wedge
\tau_U} \left( \int_{E\setminus V} (1-\phi_V(X_s)) |u|(y) N(X_s,
dy)\right) dH_s$ of $X^U$ is $\mu:= \left(\int_{E\setminus V}
(1-\phi_V(x)) |u|(x) N(x, dy)\right) d\mu_H$ and so $f=G_U \mu$.
Since by condition \eqref{e:cond1},
  $$ \mu (U)= \int_U \left(\int_{E\setminus V} (1-\phi_V (y))|u(y)| N(x, dy)\right) \mu_H(dx)
  \leq \int_U \left(\int_{E\setminus V}  |u(y)| N(x, dy)\right) \mu_H(dx) <\infty,
  $$
  we have
 $ \int_U G_U \mu (x) \mu (dx) \leq \|f \|_\infty \, \mu
  (U)<\infty$. Applying
    Lemma \ref{L:2.4} to $X^U$ yields that  $f\in \FF^U_e$. \qed

\begin{lemma}\label{L:2.5} Let $D$ be an open subset of $E$.
Every $u\in \FF_e$ that is locally bounded on $D$ satisfies
conditions \eqref{e:cond1} and \eqref{e:cond2}.
\end{lemma}

\pf  Let $u\in \FF_e$ be locally bounded on $D$. For
   any   relatively compact open sets $U, V$ with
 $\overline U \subset V \subset \overline V \subset D$, take $\phi
 \in \FF\cap C_c(D)$ such that $\phi =1$ on $U$ and $\phi =0$ on
 $V^c$. Then $u\phi \in \FF_e$ and
 \begin{eqnarray*}
  \int_{U\times  (E\setminus V)} u(y)^2 J(dx, dy)
 &= & \int_{U\times  (E\setminus V)} \left((1-\phi)u)(x)-((1-\phi) u)(y)\right)^2 J(dx,
 dy)\\
&\leq&  2\EE (u-u\phi, u-u\phi)<\infty.
\end{eqnarray*}
This together with \eqref{e:J1} implies that
$$ \int_{U\times  (E\setminus V)} |u(y)| J(dx, dy)
\leq \frac12 \int_{U\times  (E\setminus V)} \left(1+ u(y)^2\right) J(dx, dy) <\infty.
$$

Let $\phi_V\in \FF\cap C_c(D)$ be such that $0\leq \phi_V\leq 1$ with
$\phi_V=1$ on $V$. Note that $|u|\in \FF_e$ is locally bounded on $D$
and so $(1-\phi_V)|u|= |u|-\phi_V |u|\in \FF_e$. Thus it follows from \cite[Theorem 3.4.8]{CF} or \cite[Theorem 4.6.5]{FOT} that
$$ \1_U(x)\E_x \left[
\big((1-\phi_V  \big) |u| (X_{\tau_U}) \right]= \E_x \left[
\big((1-\phi_V ) |u|\big) (X_{\tau_U})\right] - (1-\phi_V) |u| \in
\FF^U_e. $$
 \qed

 \begin{lemma}\label{L:3} Let $D$ be a relatively compact open set of $E$.
  Suppose $u$ is a function in
  $\FF^D_\loc$
    that is locally bounded on $D$ and satisfies
   the condition \eqref{e:cond1}.
  Then for every $v\in C_c(D)\cap \FF$,
 the expression
 $$ \frac12 \mu^c_{\<u, v\>}(D) + \frac12 \int_{E\times E} (u(x)-u(y))(v(x)-v(y)) J(dx, dy)
+ \int_D u(x) v(x) \kappa (dx)
 $$
 is well-defined and finite; it  will still be denoted as $\EE(u, v)$.
 \end{lemma}

\pf  Clearly the first and the third terms are well defined and finite. To see that the second term
 is also well defined,
 let $U$ be a relatively compact open subset  of $D$ such that
 ${\rm supp} [v]\subset U  $. Since $u\in \FF^D_\loc$, there is $f\in \FF$ so that
 $u=f$ $m$-a.e. and hence q.e. on $U$.
  Under condition \eqref{e:cond1},
  \begin{eqnarray*}
 && \int_{E\times E} |(u(x)-u(y))(v(x)-v(y))| J(dx, dy) \\
  &\leq & \int_{U\times U} |(u(x)-u(y))(v(x)-v(y))| J(dx, dy) +
  2 \int_{U\times (E\setminus U)} |u(x)v(x)| J(dx, dy) \\
 &&     + 2 \int_U |v(x)| \int_{E\setminus U} |u(y)| J(dx, dy)  \\
  &\leq & \int_{U\times U} |(f (x) -f(y))(v(x)-v(y))| J(dx, dy)
  + 2 \| uv\|_\infty J( {\rm supp}[v], U^c) )\\
 &&          + 2   \|v\|_\infty \int_{{\rm supp}[v]\times (E\setminus U)} |u(y)| J(dx,
          dy) \\
         &<& \infty.
         \end{eqnarray*}
In the last inequality we used \eqref{e:J1} and the fact that $f, v\in \FF$.
 This proves the lemma. \qed

 \begin{thm}\label{T:4}
 Let $D$ be an open subset of $E$. Suppose that
 $u\in \FF^D_{\loc} $
  is locally bounded on $D$  satisfying
  conditions
 \eqref{e:cond1}-\eqref{e:cond2} and that
\begin{equation}\label{e:2.3}
\EE(u, v)=0 \qquad \hbox{for every } v \in C_c(D)\cap \FF .
\end{equation}
 Then $u$ is harmonic in $D$.
  If $U$ is a relatively compact open subset of $D$ so that  $\P_x(\tau_U<\infty)>0$ for q.e. $x\in U$,
 then $u(x)=\E_x \left[ u(X_{\tau_U})\right]$
 for q.e. $x\in U$.
 \end{thm}

 \pf Take  $\phi \in C_c(D)\cap \FF $
  such that $0\leq \phi \leq 1$ and $\phi =1$  in an open neighborhood $V$ of $\overline U$.
  Then $\phi u\in \FF^D$. So by \cite[Theorem 3.4.8]{CF} or \cite[Theorem 4.6.5]{FOT},
  $h_1(x):=\E_x \left[ (\phi u)(X_{\tau_U})\right]  \in \FF_e$ and
  $\phi u-h_1\in \FF_e^U$. Moreover
  \begin{equation}\label{eqn:2}
   \EE (h_1, v)=0 \qquad \hbox{for every } v\in \FF_e^U.
  \end{equation}
  Let $h_2(x):=\E_x\left[ ((1-\phi)u) (X_{\tau_U})\right]$, which is well defined
  by condition \eqref{e:cond2}. Note that by the L\'evy system of
  $X$,
  $$ f(x):= \1_U (x) \E_x\left[ \big((1-\phi)|u| \big) (X_{\tau_U})\right]
   = \1_U (x)\, \E_x \left[ \int_0^{\tau_U}\left( \int_{E\setminus V} \big((1-\phi)|u| \big) (z)
     N(X_s. dz) \right)dH_s\right].
$$
Define $\mu  (dx):= \1_D(x) \left( \int_{E\setminus V}
\big((1-\phi)|u| \big) (z) N(X_s. dz) \right) \mu_H(dx)$, which is a
smooth measure of $X^U$. In the following, for a smooth measure
$\nu$ of $X^U$, we will use $G_U\nu$ to denote $\E_x [
A^\nu_{\tau_U}]$, where $A^\nu$ is the PCAF of $X^U$ with Revuz
measure $\nu$. Using such a notation, $f=G_U \mu $.
 We claim that $\1_U h_2 \in \FF^U_e$ and for $v\in \FF^U_e$,
 \begin{equation}\label{e:2.12}
 \EE (\1_D h_2, \, v) = \int_E v(x)  \1_U(x) \left( \int_{E\setminus V}
\big((1-\phi) u \big) (z) N(X_s, dz) \right) \mu_H(dx).
\end{equation}
Define \begin{eqnarray*}
 \mu_1 (dx) &:=& \1_D(x) \left(
\int_{E\setminus V} \big((1-\phi) u^+ \big) (z) N(X_s. dz) \right)
\mu_H(dx) , \\
 \mu_2 (dx) &:=& \1_D(x) \left(
\int_{E\setminus V} \big((1-\phi) u^- \big) (z) N(X_s. dz) \right)
\mu_H(dx).
\end{eqnarray*}
Observe that
$$G_U \mu_1(x)= \E_x\left[ ((1-\phi)u^+) (X_{\tau_U})\right] \quad
  \hbox{ and } \quad G_U \mu_2(x)= \E_x\left[ ((1-\phi)u^-) (X_{\tau_U})\right]
  \qquad \hbox{for } x\in U.
  $$
Clearly $G_U \mu_1 \leq G_U \mu$. For $j\geq 1$, let $F_j:=\{x\in U:
G_U \mu_1 (x) \leq j\}$, which is a finely closed subset of $U$.
Define $\nu_j:=\1_{F_j} \mu_1$. Then for $x\in F_j$, $G_U \nu_j (x)
\leq G_U \mu_1 (x) \leq j$, which for $x\in U\setminus F_j$,
$$ G_U \nu_j(x) = \E_x  \left[ G_U \nu_j( X_{\sigma_{F_j}}) \right] \leq j.
$$
In other words, we have $G_U \nu_j \leq j\wedge G_U \mu_1 \leq j
\wedge f$. As both $G_U \nu_j$ and $j\wedge f$ are excessive
functions of $X^U$ and $m(U)<\infty$, we have by \cite[Theorem 1.1.5
and Lemma 1.2.3]{CF} that $\{G_U \nu_j, \ j\wedge G_U \mu\}\subset
\FF^U$ and
$$ \EE (G_U \nu_j, \, G_U \nu_j) \leq \EE (j\wedge f, \,
j\wedge f) \leq \EE (f, \, f ) <\infty.
$$
Moreover,   for each $j\geq 1$, we have by \cite[Theorem 4.1.1]{CF}
or \cite[Theorem 5.1.3]{FOT} that
\begin{eqnarray*} \EE
(G_U \nu_j, \, G_U \nu_j)
  &=& \lim_{t\to 0} \frac1t \int_E G_U (\nu_j (x) -P^U_t G_U \nu_j(x))
 G_U \nu_j(x) m(dx) \\
 &=& \lim_{t\to 0} \frac1t \int_E \E_x \left[ A^{\nu_j}_{t\wedge \tau_U} \right]
 G_U \nu_j(x) m(dx) \\
 &=& \int_U G_U \nu_j(x) \, \1_{F_j}(x) \mu_1 (dx),
\end{eqnarray*}
which increases to $\int_U G_U \mu_1 (x) \mu_1 (dx)$. Consequently,
$\int_U G_U \mu_1 (x) \mu_1 (dx)\leq \EE (f, f)<\infty$. So we have by
Lemma \ref{L:2.4} applied to $X^U$ that $G_U \mu_1 \in \FF^U_e$ with
$\EE (G_U \mu_1, v)=\int_U v(x) \mu_1(dx)$ for every $v\in \FF^U_e$.
Similarly we have $G_U \mu_2\in \FF^U_e$ with $\EE (G_U \mu_2,
v)=\int_U v(x) \mu_2(dx)$ for every $v\in \FF^U_e$. It follows that
$\1_U h_2= G_U \mu_1 -G_U\mu_2\in \FF^U_e$ and claim \eqref{e:2.12}
is established.

 As $h_2=\1_U h_2 +(1-\phi)u$ and $(1-\phi)u$ satisfies condition
 \eqref{e:cond1}, we have by Lemma \ref{L:3} and \eqref{e:2.12}
 that for every $v\in C_c(U)\cap \FF$,
 \begin{eqnarray}
 \EE (h_2, v)&=& \EE (1_U h_2, v)+\EE ((1-\phi)u, v) \nonumber \\
 &=& \int_{E\times
E} v(x) (1-\phi (y)) u(y) N(x, dy) \mu_H (dx) - \int_{E\times E}
v(x) (1-\phi (y)) u(y) N(x, dy) \mu_H (dx)
 \nonumber \\
 &=& 0  . \label{eqn:3}
 \end{eqnarray}
 This combining with \eqref{eqn:2} and condition \eqref{e:2.3} proves that
\begin{equation}\label{eqn:4}
\EE(u-h_1-h_2, v)=0 \qquad \hbox{for every }  v\in C_c(U)\cap
 \FF.
 \end{equation}
 Since  $u-(h_1+h_2)=(\phi u -h_1)- \1_D h_2  \in \FF_e^U$ and $C_c(U)\cap \FF$
  is $\EE$-sense in $\FF^U_e$,  the above display holds for every $v\in \FF^U_e$.
  In particular, we have
  \begin{equation}\label{e:2.8}
   \EE (u-h_1-h_2, \, u-h_1-h_2)=0.
  \end{equation}
  By Lemma \ref{L:2},
  $u(X_t)-h_1(X_t)-h_2(X_t)$ is a bounded $\P_x$-martingale
  for q.e. $x\in E$.
  As
  $$h_1(x)+h_2(x)=\E_x \left[ u(X_{\tau_U})\right] \qquad \hbox{for } x\in U,
  $$
   the above implies that
  $t\mapsto  u(X_{t\wedge \tau_U})$ is a uniformly integrable $\P_x$-martingale
  for q.e. $x\in U$.
  If $\P_x(\tau_U <\infty )>0$ for q.e. $x\in U$,
  applying Lemma \ref{L:2} to the Dirichlet form $(\EE , \FF^U)$, we have
   $u-h_1-h_2=0$ q.e. on $U$
  and so
  $u(x)=\E_x \left[ u(X_{\tau_U})\right]$ for q.e. $x\in U$. This completes the proof of
   the theorem. \qed

\begin{remark}\label{R:5} \rm  \begin{description} \item{(i)}
The principal difficulty in above proof is establishing
\eqref{eqn:4} and that $u-(h_1+h_2)\in \FF^U_e$ for general $u\in
\FF^D_\loc$ satisfying conditions \eqref{e:cond1} and
\eqref{e:cond2}. If $u$ is assumed a priori to be in $\FF_e$, these
facts and therefore the theorem itself are then much easier to
establish. Note that when $u\in \FF_e$, it follows immediately from
\cite[Theorem 3.4.8]{CF} or \cite[Theorem 4.6.5]{FOT} that $h_1+h_2=
\E_x [u(X_{\tau_U})]\in \FF_e$ enjoys property \eqref{eqn:4} and
$u-(h_1+h_2)\in \FF^U_e$. Therefore \eqref{e:2.8} holds and
consequently $u$ is harmonic in $D$.

\item{(ii)} If we assume that the process $X$ (or equivalently $(\EE, \FF)$)
is $m$-irreducible and that $U^c$ is not $m$-polar, then    $\P_x
(\tau_U<\infty )>0$ for q.e. $x\in U$ (cf. \cite[Theorem 3.5.6]{CF}
or \cite{FOT}).
\end{description}
\end{remark}

 \begin{thm}\label{T:6}
  Suppose $D$ is an open set of $E$ with $m(D)<\infty$
  and $u$ is  a function on $E$
  satisfying the  condition \eqref{e:cond1}
  so that $u\in L^\infty  (D; m)$ and
  that $\{u(X_{t\wedge \tau_D}), t\geq 0\}$ is a uniformly integrable
  $\P_x$-martingale for q.e. $x\in E$.  Then
  \begin{equation}\label{e:1.5}
  u\in \FF^D_{\loc} \qquad \hbox{and} \qquad
  \EE(u, v)=0 \quad \hbox{for every } v \in C_c(D)\cap \FF .
  \end{equation}
 \end{thm}

  \pf  As  for q.e. $x\in E$, $\{u(X_{t\wedge \tau_D}), t\geq 0\}$ is a
  uniformly integrable  $\P_x$-martingale,
    $ u(X_{t\wedge \tau_D})$  converges
        in $L^1(\P_x)$ as well as $\P_x$-a.s. to some random variable
        $\xi$. By considering $\xi^+$, $\xi^-$ and $u_+:=\E_x [\xi^+]$,
        $u_-(x):=\E_x [\xi^-]$ separately,  we may and do assume without loss
        of generality that $u\geq 0$. Note that $\xi \1_{\{\tau_D<\infty\}}
         = u(X_{\tau_D})$. Define $u_1(x):=\E_x \left[ u(X_{\tau_D})\right]$
         and $u_2(x):=\E_x [ \xi \1_{\{\tau_D = \infty\}}] =u-u_1$.

  Let $\{P^D_t, t\geq 0\}$ denote the transition semigroup of the subprocess $X^D$.
   Then for q.e. $x\in D$ and every $t>0$,
    by the Markov property of $X^D$,
$$  P^D_t u_2(x)=\E_x \left[ u_2(X_t), t<\tau_D\right] =
  \E_x \left[ \xi \1_{\{\tau_D=\infty\}} \cdot \theta_t, t<\tau_D\right] =u_2(x) .
  $$
  Since $u_2 \in L^2(D; m)$, by \eqref{e:1.1}-\eqref{e:1.2}
  \begin{equation}\label{e:u2}
  u_2\in \FF^D \quad \hbox{with} \quad  \EE (u_2, u_2)=0.
  \end{equation}
   On the other hand,
  $$ P^D_t u(x)=\E_x \left[ u(X_t), t<\tau_D\right] =
  \E_x \left[ u(X_{\tau_D}), t<\tau_D\right] \leq u(x) .
  $$
  Let $\{D_n, n\geq 1\}$ be an increasing sequence of relatively compact open subsets
  of $D$ with $\cup_{n\geq 1} D_n=D$ and define
  $$ \sigma_n:=\inf \left\{t\geq 0: X^D_t\in D_n\right\}.
  $$
   Let $e_n(x)=\E_x \left[ e^{-\sigma_n}\right]$,
  $x\in D$, be the 1-equilibrium potential of $D_n$ with respect to the subprocess $X^D$.
  Clearly $e_n \in \FF^D$ is 1-excessive with respect to the process $X^D$,
  $e_n(x)=1$ q.e. on $D_n$. Let $a:=\| \1_D u\|_\infty$.
  Then  for every $t>0$,
  $$ e^{-t} P^D_t( (a e_n)\wedge u)(x) \leq ((a e_n) \wedge u) (x) \qquad \hbox{for q.e. } x\in D.
  $$
By \cite[Lemma 1.2.3]{CF} or \cite[Lemma 8.7]{S}, we have
$(ae_n)\wedge u\in \FF^D$ for every $n\geq 1$. Since $(a e_n)\wedge
u=u$ $m$-a.e. on $D_n$, we have  $u\in \FF^D_{\loc} $.

Let $U$ be a relatively
compact open subset of $D$. Let $\phi \in C_c(D)\cap \FF$ so that $0\leq \phi \leq 1$
and  $\phi =1$ in an open neighborhood $V$ of $\overline U$.
Define for $x\in E$,
$$ h_1(x):= \E_x \left[ (\phi u)(X_{\tau_U})\right]
\quad \hbox{and} \quad  h_2(x):= \E_x \left[ ((1-\phi)
u)(X_{\tau_U})\right].
$$
Then $u_1=h_1+h_2$ on $E$. Since $\phi u\in \FF$, we know as in
(\ref{eqn:2}) that $h_1\in \FF_e$ and
$$ \EE (h_1, v)=0 \qquad \hbox{for every } v\in \FF_e^U.
$$
By the same argument as that for (\ref{eqn:3}), we have
$$ \EE (h_2, v)=0 \qquad \hbox{for every } v\in \FF_e^U.
$$
These together with \eqref{e:u2} in particular implies that
$$ \EE (u, v)=\EE(h_1+h_2+u_2, v)=0 \qquad \hbox{for every } v\in C_c(U)\cap \FF .
$$
Since $U$ is an arbitrary relatively compact subset of $D$, we have
$$ \EE (u, v)= 0 \qquad \hbox{for every } v\in C_c(D)\cap \FF .
$$
This completes the proof. \qed

\begin{remark}\label{R:2.7} \rm As mentioned in the Introduction, the
principal difficulty for the proof of the above theorem is
establishing that a   function $u$ harmonic in $D$ is in
$\FF^D_\loc$ with $\EE(u, v)=0$ for every $v\in \FF \cap C_c(D)$. If
a priori $u$ is assumed to be in $\FF_e$, then Theorem \ref{T:6} is
easy to establish. In this case, it follows from \cite[Theorem
3.4.8]{CF} or \cite[Theorem 4.6.5]{FOT} that $u_1=h_1+h_2 = \E_x [
u(X_{\tau_U})]\in \FF_e$ and that the second property of
\eqref{eqn:4} holds. This together with \eqref{e:u2} immediately
implies that $u$ enjoys \eqref{e:1.5}. (See also Proposition 2.5 of
\cite{BBKT} for this simple case but under an additional assumption that
$1\in \FF$ with $\EE(1, 1)=0$.)
\end{remark}

\medskip

Combining Theorems \ref{T:4} and \ref{T:6},
we have the following.

\begin{thm}\label{T:7} Let  $D$ be an open subset of $E$.
Suppose that $u$ is a function on $E$ that is locally bounded on $D$
and satisfies conditions \eqref{e:cond1} and \eqref{e:cond2}. Then
\begin{description}
\item{\rm (i)} $u$ is harmonic in $D$ if and only if condition \eqref{e:1.5} holds.
 \item{\rm (ii)} Assume  that for every relatively compact open subset $U$ of $D$,
 $\P_x (\tau_U < \infty)>0$ for q.e. $x\in U$.
 {\rm (}By Remark \ref{R:5}(ii), this condition is satisfied if $(\EE, \FF)$ is
 $m$-irreducible.{\rm )} Then $u$ is harmonic in $D$  if and only if
 for every relatively compact subset $U$ of $D$, $u(X_{\tau_U})\in
 L^1 (\P_x)$ and $u(x)= \E_x \left[ u(X_{\tau_U})\right]$ for
 q.e. $x\in U$.
 \end{description}
\end{thm}

\medskip

\begin{example}\label{E:8} \rm   (Stable-like process on $\R^d$)
Consider the following   Dirichlet form $(\EE, \FF)$ on $L^2(\R^d,
dx)$, where
\begin{eqnarray*}
 \FF &=& W^{\alpha/2, 2}(\R^d):=\left\{ u\in L^2(\R^d; dx): \ \int_{\R^d\times \R^d}
 (u(x)-u(y))^2
 \frac{1}{|x-y|^{d+\alpha}} \, dxdy<\infty \right\} , \\
  \EE(u, v)&=& \frac12 \int_{\R^d\times \R^d} (u(x)-u(y))(v(x)-v(y))
\frac{c(x, y)}{|x-y|^{d+\alpha}}\, dx dy \qquad \hbox{for } u, v\in
\FF.
\end{eqnarray*}
Here $d\geq 1$, $\alpha \in (0, 2)$ and $c(x, y)$ is a symmetric
function in $(x, y)$ that is bounded between two positive constants.
In literature, $W^{\alpha, 2}(\R^d)$ is called the Sobolev space on
$\R^d$ of fractional order $(\alpha/2, 2)$. For an open set
$D\subset \R^d$, $W^{\alpha, 2}(D)$ is similarly defined as above
but with $D$ in place of $\R^d$. It is easy to check that $(\EE,
\FF)$ is a regular Dirichlet form on $L^2(\R^d; dx)$ and its
associated symmetric Hunt process $X$ is called symmetric
$\alpha$-stable-like process on $\R^d$, which is
 studied in \cite{CK}. The process $X$ has strictly positive
jointly continuous transition density function $p(t, x, y)$ and
hence is irreducible. Moreover, there is constant $c>0$ such that
\begin{equation}\label{e:2.18}
p(t, x, y) \leq c \, t^{-d/\alpha} \qquad \hbox{for } t>0 \hbox{ and
} x, y \in \R^d
\end{equation}
and consequently by \cite[Theorem 1]{Ch},
\begin{equation}\label{e:2.19}
\sup_{x\in U} \E_x [ \tau_U] <\infty
\end{equation}
for any open set $U$ having finite Lebesgue measure.
 When $c(x, y)$ is constant, the process $X$ is nothing but
the rotationally symmetric $\alpha$-stable process on $\R^d$. In
this example, the jumping measure
$$J(dx, dy)= \frac{c(x, y)}{|x-y|^{d+\alpha}} \, dx dy.
$$
Hence for any non-empty open set $D\subset \R^d$,  condition \eqref{e:cond1} is satisfied if and only if
$(1\wedge |x|^{-d-\alpha} ) u(x) \in L^1  (\R^d)$. Moreover, for
such a function $u$ and relatively compact open sets $U, V$ with
$\overline U\subset V \subset  \overline V \subset D$, by L\'evy
system of $X$,
\begin{eqnarray}\label{e:2.20}
\sup_{x\in U} \E_x \left[ (\1_{V^c} |u|) (X_{\tau_U})\right] &=&
\sup_{x\in U} \E_x \left[ \int_0^{\tau_U} \left( \int_{V^c}  \frac{
c(X_s, y) \, |u(X_s)| }{|X_s-y|^{d+\alpha}} dy \right) ds \right]
\nonumber \\
&\leq & \left( c\, \int_{\R^d} (1\wedge |y|^{-d-\alpha} ) |u(y)|dy
\right) \sup_{x\in U} \E_x [ \tau_U] < \infty .
\end{eqnarray}
In other words, for this example, condition \eqref{e:cond3} and hence \eqref{e:cond2} is a
consequence of \eqref{e:cond1}.
 So Theorem \ref{T:6} says that for an open set $D$ and
 a function $u$ on $\R^d$ that is locally bounded on $D$ with
 $(1\wedge |x|^{-d-\alpha} ) u(x) \in L^1  (\R^d)$, the following are equivalent.

 \begin{description}
 \item{(i)} $u$ is harmonic in $D$;
 \item{(ii)} For every relatively compact subset $U$ of $D$, $u(X_{\tau_U})\in
 L^1 (\P_x)$ and $u(x)= \E_x \left[ u(X_{\tau_U})\right]$ for
 q.e. $x\in U$;
 \item{(iii)}  $u\in \FF^D_\loc =W^{\alpha, 2}_\loc (D)$
  and
 $$ \int_{\R^d\times \R^d} (u(x)-u(y))(v(x)-v(y)) \frac{c(x, y)}{|x-y|^{d+\alpha}}
  \, dxdy =0 \qquad \hbox{for every } v\in C_c(D)\cap
  W^{\alpha/2, 2}(\R^d).
  $$
  \end{description}
  \qed
\end{example}

\begin{example}\label{E:9} \rm (Diffusion process on a locally compact separable metric space)\
 Let $(\EE, \FF)$ be a local regular Dirichlet form on $L^2(E; m)$, where $E$
 is a locally compact separable metric space, and $X$ is its associated Hunt process.
 In this case, $X$ has continuous sample paths and so
 the jumping measure $J$ is null (cf. \cite{FOT}). Hence  conditions
 \eqref{e:cond1} and \eqref{e:cond2}
 are automatically satisfied. Let $D$ be an open subset of
  $E$
 and $u$ be a function
 on $E$ that is locally bounded in $D$. Then by Theorem \ref{T:7},
 $u$ is harmonic in $D$
 if and only if condition \eqref{e:1.5} holds.

  Now consider the
 following special case: $E=\R^d$ with $d\geq 1$, $m(dx)$ is the Lebesgue measure
 $dx$ on $\R^d$, $\FF=W^{1,2}(\R^d):=\left\{u\in L^2(\R^d; dx) \mid
 \nabla u \in L^2(\R^d; dx) \right\}$ and
 $$ \EE (u, v) = \frac12 \sum_{i,j=1}^d \int_{\R^d}  a_{ij}(x)
 \frac{\partial u(x)}{\partial x_i}  \frac{\partial v (x)}{\partial
 x_j} \, dx  \qquad \hbox{for } u, v \in W^{1,2}(\R^d),
 $$
 where $(a_{ij}(x))_{1\leq i, j\leq d}$ is a $d\times d$-matrix
 valued measurable function on $\R^d$ that is uniformly elliptic
 and bounded.
In literature, $W^{1, 2}(\R^d)$ is  the Sobolev space on $\R^d$ of
order $(1, 2)$. For an open set $D\subset \R^d$, $W^{1, 2}(D)$ is
similarly defined as above but with $D$ in place of $\R^d$.
 Then $(\EE, \FF)$ is a regular local Dirichlet form on $L^2(\R^d; dx)$
 and its associated Hunt process $X$ is a conservative diffusion on
 $\R^d$ having jointly continuous transition density function.
 Let $D$ be an open set in $\R^d$. Then by Theorem \ref{T:7}, the following are equivalent
 for a locally bounded
 function $u$ on $D$.
 \begin{description}
 \item{(i)} $u$  is harmonic in $D$;
 \item{(ii)} For every relatively compact open subset $U$ of
 $D$, $u(X_{\tau_U})\in L^1(\P_x)$ and $u(x)=\E_x \left[
  u(X_{\tau_U})
  \right]$ for q.e. $x\in U$;
 \item{(iii)} $u\in W^{1,2}_\loc (D)$ and
 $\displaystyle \sum_{i,j=1}^d \int_{\R^d}   a_{ij}(x)
 \frac{\partial u(x)}{\partial x_i}  \frac{\partial v (x)}{\partial
 x_j} \, dx =0$  \ for every  $ v \in C_c(D)\cap W^{1,2}(\R^d)$.
\end{description}
In fact, in this case, it can be shown that every (locally bounded)
harmonic function has a continuous version.  \qed
 \end{example}

\begin{example}\label{E:10} \rm (Diffusions with jumps on $\R^d$)
Consider the following Dirichlet form $(\EE, \FF)$, where
$\FF=  W^{1,2}(\R^d)$ and
\begin{eqnarray*}
  \EE(u, v)&=& \frac12 \sum_{i,j=1}^d \int_{\R^d}   a_{ij}(x)
 \frac{\partial u(x)}{\partial x_i}  \frac{\partial v (x)}{\partial
 x_j} \, dx \\&& + \frac12 \int_{\R^d\times \R^d} (u(x)-u(y))(v(x)-v(y))
\frac{c(x, y)}{|x-y|^{d+\alpha}}\, dx dy \qquad \hbox{for } u, v\in
 W^{1, 2}(\R^d).
\end{eqnarray*}
Here $d\geq 1$, $(a_{ij}(x))_{1\leq i, j\leq d}$ is a $d\times d$-matrix
 valued measurable function on $\R^d$ that is uniformly elliptic
and bounded, $\alpha \in (0, 2)$ and $c(x, y)$ is a symmetric
function in $(x, y)$ that is bounded between two positive constants.
It is easy to check that $(\EE, \FF)$ is a regular Dirichlet form on
$L^2(\R^d; dx)$. Its associated symmetric Hunt process $X$ has both
the diffusion and jumping components. Such a process has recently
been studied in \cite{CK2}. It is shown there the process $X$ has
strictly positive jointly continuous transition density function
$p(t, x, y)$ and hence is irreducible. Moreover, a sharp two-sided
estimate is obtained in \cite{CK2} for $p(t, x, y)$. In particular,
there is a constant $c>0$ such that
$$ p(t, x, y) \leq c \left( t^{-d/\alpha} \wedge t^{-d/2}\right)
\qquad \hbox{for } t>0 \hbox{ and } x, y\in \R^d.
$$
Note that when $(a_{ij})_{1\leq i, j\leq d}$ is the identity matrix
and $c(x, y)$ is constant, the process $X$ is nothing but the
symmetric L\'evy process that is the independent sum of a Brownian
motion and a rotationally symmetric $\alpha$-stable process on
$\R^d$. In this example, the jumping measure
$$J(dx, dy)= \frac{c(x, y)}{|x-y|^{d+\alpha}} \, dx dy.
$$
Hence for any non-empty open set $D\subset \R^d$,
condition \eqref{e:cond1} is satisfied if and only if
$(1\wedge |x|^{-d-\alpha} ) u(x) \in L^1  (\R^d)$. By the same
reasoning as that for \eqref{e:2.20}, we see that for this example,
condition \eqref{e:cond3} and hence \eqref{e:cond2} is implied by condition \eqref{e:cond1}.
So Theorem \ref{T:6} says that for an open set $D$ and
 a function $u$ on $\R^d$ that is locally bounded on $D$ with
 $(1\wedge |x|^{-d-\alpha} ) u(x) \in L^1  (\R^d)$, the following are equivalent.

 \begin{description}
 \item{(i)} $u$ is harmonic in $D$ with respect to $X$;
 \item{(ii)} For every relatively compact subset $U$ of $D$, $u(X_{\tau_U})\in
 L^1 (\P_x)$ and $u(x)= \E_x \left[ u(X_{\tau_U})\right]$ for
 q.e. $x\in U$;
 \item{(iii)} $u\in W^{1,2}_\loc (D)$ such that   for every $ v\in C_c(D)\cap W^{1,2}(\R^d)$,
 $$  \sum_{i,j=1}^d  \int_{\R^d}  a_{ij}(x)
 \frac{\partial u(x)}{\partial x_i}  \frac{\partial v (x)}{\partial
 x_j} \, dx+ \int_{\R^d\times \R^d} (u(x)-u(y))(v(x)-v(y)) \frac{c(x, y)}{|x-y|^{d+\alpha}}
  \, dxdy =0 .
  $$
  \end{description}
\qed
\end{example}

\begin{remark}\label{R:11} \rm It is possible to extend the results of this paper
to   a general $m$-symmetric right process $X$
 on a Lusin space, where  $m$ is a positive $\sigma$-finite
measure  with full topological support on $E$. In this case, the
Dirichlet $(\EE, \FF)$ of $X$ is   a quasi-regular Dirichlet form on
$L^2(E; m)$. By \cite{CMR}, $(\EE, \FF)$ is quasi-homeomorphic to a
regular Dirichlet form on a locally compact separable metric space.
So the results of this paper can be extended to the quasi-regular
Dirichlet form setting, by using this quasi-homeomorphism. However
since the notion of open set is not invariant under
quasi-homeomorphism, some modifications are   needed. We need to
replace open set $D$ in Definition \ref{D:1} by quasi-open set $D$.
Similar modifications are needed for conditions \eqref{e:cond1} and
\eqref{e:cond2} as well. We say a function $u$ is harmonic in a
quasi-open set $D\subset E$ if for every quasi-open subset $U\subset
D$ with $\overline U \cap F_k \subset D$ for every $k\geq 1$, where
$\{F_k, k\geq 1\}$ is an $\EE$-nest consisting of compact sets,
$t\mapsto u(X_{t\wedge \tau_{U\cap F_k}})$ is a uniformly integrable
$\P_x$-martingale for q.e. $x\in U\cap F_k$ and for every $k\geq 1$.
 The
local Dirichlet space $\FF_\loc^D$ needs to be replaced by
\begin{eqnarray*}
  \stackrel{ \circ } {{\FF}_{\loc}^D} &=& \Big\{u: \hbox{ there is an increasing
       sequence of quasi--open sets }  \{D_n\}
       \hbox{ with } \bigcup^\infty_{n=1} D_n = D
        \hbox{ q.e.} \\
     &&\hskip 0.4truein \hbox{and a
         sequence }  \{u_n\} \subset \FF^D
         \hbox{ such that }
        u = u_n\  m \hbox{-a.e. on } \ D_n \Big\}.
\end{eqnarray*}
Condition \eqref{e:1.5} should be replaced by
\begin{equation}\label{e:1.7}
 u\in \stackrel{ \circ } {{\FF}_{\loc}^D} \qquad \hbox{and} \qquad
\EE (u, v)=0 \quad \hbox{for every } v\in \FF \hbox{ with } \EE \hbox{-supp[$v$]}
\subset D.
\end{equation}
Here $\EE$-supp[$u$] is the smallest quasi-closed set that $u$
vanishes $m$-a.e. on its complement. We leave the details to
interested readers. \qed
\end{remark}

\noindent {\bf Acknowledgement.} The author thanks Rich Bass and
Takashi Kumagai for helpful discussions. He also thanks Rongchan Zhu
for helpful comments.

\vskip 0.1truein

\small

\begin{singlespace}

\end{singlespace}

\bigskip

Department of Mathematics, University of Washington, Seattle, WA
98195, USA.

   Email:  {\texttt zchen@math.washington.edu}

\end{doublespace}

\end{document}